 \documentclass[12pt]{article}
 \usepackage{indentfirst, latexsym,bm}
 \usepackage{amsmath}
 \usepackage{pifont}
 \usepackage{amsfonts}
 \usepackage{mathrsfs}
 \usepackage{array}
 \usepackage{multirow}
 \usepackage{graphicx}
 \usepackage{subfigure}
 \usepackage{picinpar}

 \usepackage{booktabs}
 \usepackage{threeparttable}
 \usepackage{mathrsfs,amsfonts,amsmath}
 \usepackage[dvips]{color}

\newcommand\email[1]{\href{mailto:#1}{ \nolinkurl{#1}}}

 \newtheorem{theorem}{Theorem}[section]
 \newtheorem{definition}[theorem]{Definition}
 \newtheorem{lemma}[theorem]{Lemma}
 \newtheorem{corollary}[theorem]{Corollary}
 \newtheorem{proposition}[theorem]{Proposition}
 \newtheorem{remark}[theorem]{Remark}
 \newtheorem{condition}[theorem]{Condition}
 \newtheorem{example}{Example}[section]

 \def\blemma{\begin{lemma}}\def\elemma{\end{lemma}}
 \def\bproposition{\begin{proposition}}\def\eproposition{\end{proposition}}
 \def\ttheorem{\begin{theorem}}\def\etheorem{\end{theorem}}
 \def\bcorollary{\begin{corollary}}\def\ecorollary{\end{corollary}}
 \def\bremark{\begin{remark}}\def\eremark{\end{remark}}
 \def\bcondition{\begin{condition}}\def\econdition{\end{condition}}

 \def\benumerate{\begin{enumerate}}\def\eenumerate{\end{enumerate}}
 \def\bitemize{\begin{itemize}}\def\eitemize{\end{itemize}}
 \def\itm{\item}

 \def\beqlb{\begin{eqnarray}}\def\eeqlb{\end{eqnarray}}
 \def\beqnn{\begin{eqnarray*}}\def\eeqnn{\end{eqnarray*}}
 \def\ar{\!\!\!&}

 \def\mcr{\mathscr}\def\mbb{\mathbb}\def\mbf{\mathbf}

 \def\proof{\noindent{\it Proof.~~}}\def\qed{\hfill$\Box$\medskip}

\begin{document}

\noindent{(Version: 2016/01/19)}

\bigskip

\noindent{\LARGE\bf Continuous-state branching processes}

\noindent{\LARGE\bf in L\'evy random environments\footnote{\,Supported by NSFC
(No.\,11131003, No.\,~11371061 and No.\,11531001).}}

\bigskip

\noindent{Hui He, Zenghu Li and Wei Xu\footnote{\,Corresponding author.}}

\smallskip

\noindent{School of Mathematical Sciences, Beijing Normal University,}

\noindent{Beijing 100875, People's Republic of China}

\noindent{E-mails: \tt hehui@bnu.edu.cn, lizh@bnu.edu.cn, xuwei@mail.bnu.edu.cn}

\bigskip

\textit{Abstract.} A general continuous-state branching processes in random environment (CBRE-process) is defined as the strong solution of a stochastic integral equation. The environment is determined by a L\'evy process with no jump less than $-1$. We give characterizations of the quenched and annealed transition semigroups of the process in terms of a backward stochastic integral equation driven by another L\'evy process determined by the environment. The process hits zero with strictly positive probability if and only if its branching mechanism satisfies Grey's condition. In that case, a characterization of the extinction probability is given using a random differential equation with singular terminal condition. The strong Feller property of the CBRE-process is established by a coupling method. We also prove a necessary and sufficient condition for the ergodicity of the subcricital CBRE-process with immigration.

\medskip


\medskip



\section{Introduction}

\setcounter{equation}{0}

Galton-Watson processes in random environments (GWRE-processes) were introduced by Smith (1968) and Smith and Wilkinson (1969) as extensions of classical Galton-Watson processes (GW-processes). Those extensions possess many interesting new properties such as phase transitions of survival probabilities. For instance, that different regimes for the survival probability arise in the subcritical regime. For recent results on the speed of decay of the survival probability, the reader may refer to Afanasy'ev et al.\ (2005, 2012), Guivarc¡¯h and Liu (2001) and Vatutin (2004) and the references therein.

Let $\{W(t): t\ge 0\}$ be a Brownian motion and $\{S(t): t\ge 0\}$ a Brownian motion with drift. We assume the two processes are independent. By the Yamada-Watanabe theorem, for any constants $c\ge 0$, $\sigma\ge 0$ and $b\in \mbb{R}$ there is a unique positive solution to the stochastic differential equation:
 \beqlb\label{eq1.1}
dX(t) = \sqrt{2cX(t)}dW(t) - bX(t)dt + \sigma X(t)dS(t).
 \eeqlb
The solution $\{X(t): t\ge 0\}$ is called a continuous-state branching diffusion in random environment (CBRE-diffusion). The environment here is determined by the process $\{S(t): t\ge 0\}$. It was proved in Kurtz (1978) that the CBRE-diffusion arises as the limit of a sequence of suitably rescaled GWRE-processes; see also Helland (1981). A diffusion approximation of the GWRE-process was actually conjectured by Keiding (1975). It turns out that the CBRE-diffusion is technically more tractable than the GWRE-process. In the particular case of $\sigma=0$, the CBRE-diffusion reduces to the well-known Feller branching diffusion, which belongs to an important class of positive Markov processes called continuous-state branching processes (CB-processes); see Feller (1951), Ji\v{r}ina (1958) and Lamperti (1967a, 1967b).

In the work of B\"oinghoff and Hutzenthaler (2012), it was shown that the survival probability of the CBRE-diffusion can be represented explicitly in terms of an exponential functional of the environment process $\{S(t): t\ge 0\}$. Based on the representation, B\"oinghoff and Hutzenthaler (2012) gave an exact characterization for the decay rate of the survival probability of the CBRE-diffusion in the critical and subcritical cases. The results of B\"oinghoff and Hutzenthaler (2012) are more complete than the corresponding results for the GWRE-processes in the sense that they calculated the accurate limiting constants. In addition, they characterized the CBRE-process conditioned to never go extinct and established a backbone construction for the conditioned process. See also Hutzenthaler (2011) for some related results.

Continuous-state branching processes with immigration (CBI-processes), which generalize the CB-processes, were introduced by Kawazu and Watanabe (1971) as rescaling limits of Galton-Watson processes with immigration (GWI-processes); see also Aliev and Shchurenkov (1982) and Li (2006, 2011). Let $b\in \mbb{R}$ and $c\ge 0$ be given constants. Let $m(dz)$ be a Radon measure on $(0,\infty)$ satisfying $\int_0^\infty (z\land z^2)m(dz)< \infty$. Suppose that $\{W(t): t\ge 0\}$ is a Brownian motion, $\{\eta(t): t\ge 0\}$ is an increasing L\'evy process with $\eta(0)=0$ and $\tilde{M}(ds,dz,du)$ is a compensated Poisson random measure on $(0,\infty)^3$ with intensity $dsm(dz)du$. We assume those three noses are independent of each other. By Theorems~5.1 and~5.2 in Dawson and Li (2006), there is a unique positive strong solution to
 \beqlb\label{eq1.2}
X(t) \ar=\ar X(0) - b\int_0^t X(s)ds + \int_0^t\sqrt{2cX(s)}dW(s) \cr
 \ar\ar
+ \int_0^t\int_0^{\infty}\int_0^{X(s-)} z \tilde{M}(ds,dz,du) + \eta(t).
 \eeqlb
It was shown in Dawson and Li (2006) that the solution $\{X(t): t\ge 0\}$ is a CBI-processes; see also Fu and Li (2010). The process $\{\eta(t): t\ge 0\}$ describes the inputs of the immigrants. Here and in the sequel, we understand $\int_a^b = \int_{(a,b]}$ and $\int_a^\infty = \int_{(a,\infty)}$ for any $a\le b\in \mbb{R}$.

A class of continuous-state branching processes with catastrophe was introduced by Bansaye et al.\ (2013), who pointed out that the processes can be identified as continuous-state branching processes in random environment (CBRE-processes) with the environment given by a L\'evy process with bounded variation. Those authors gave a criticality classification of their CBRE-processes according to the long time behavior of the environmental L\'evy process. They also characterized the Laplace exponent of the processes using a backward ordinary differential equation involving the environment process. For stable branching CBRE-processes, Bansaye et al.\ (2013) calculated explicitly the survival probability and characterized its decay rate in the critical and subcritical cases. In addition, they showed some interesting applications of their results to a cell infection model. The results of Bansaye et al.\ (2013) were extended in Palau and Pardo (2015a) to the case where the environment was given by a Brownian motion with drift.

The CBRE-processes studied in Bansaye et al.\ (2013) can be generalized to continuous-state branching processes with immigration in random environment (CBIRE-processes). Let $\{L(t): t\ge 0\}$ be a L\'evy process with no jump less than $-1$ and assume it is independent of the three noises in (\ref{eq1.2}). It is natural to define a CBIRE-process $\{Y(t): t\ge 0\}$ by the stochastic equation
 \beqlb\label{eq1.3}
Y(t) \ar=\ar Y(0) - b\int_0^t Y(s)ds + \int_0^t\sqrt{2cY(s)}dW(s) + \eta(t) \cr
 \ar\ar
+ \int_0^t\int_0^{\infty}\int_0^{Y(s-)} z \tilde{M}(ds,dz,du) + \int_0^t Y(s-)dL(s).
 \eeqlb
Here the impact of the environment is represented by the L\'evy process $\{L(t): t\ge 0\}$. The existence and uniqueness of the positive strong solution to the above equation follow from the result of Dawson and Li (2012); see Section~2 for the details. In the very recent work of Palau and Pardo (2015b), a further generalization of the CBIRE-process was introduced by considering a competition mechanism. The authors gave a direct construction of their model by arguments similar to those in Dawson and Li (2012) and proved some results on the long time behavior of the process.

The purpose of this paper is to study the basic structures of the CBRE- and CBIRE-processes. In Section~2, we introduce two random cumulant semigroups, which are important tools in the study of those processes. The semigroups are defined in terms of a backward stochastic equations driven by a L\'evy process. The existence of them follows from a general result in Li (2011) on Dawson-Watanabe superprocesses. In Section~3, we construct the CBRE-process by applying a theorem in Dawson and Li (2012) on stochastic equations driven by time-space noises. Then we give characterizations of the quenched and annealed transition probabilities of the CBRE-process. In Section~4, we show the CBRE-process hits zero with strictly positive probability if and only if its branching mechanism satisfies Grey's condition. In that case, we give a characterization of the extinction probabilities by the solution of a random differential equation with singular terminal condition. The strong Feller property of the CBRE-process is established by a coupling method. Some of the results are extended in Section~5 to CBIRE-processes. In addition, we give a necessary and sufficient condition for the ergodicity of subcritical CBIRE-processes. Most of the results here are obtained or presented using stochastic equations driven by L\'evy processes, which are more elegant than those in the classical discrete setting.

\medskip

\textbf{Acknowledgements.} We thank Professors Yueyun Hu and Zhan Shi for helpful discussions on branching processes in random environments.


\section{Random cumulant semigroups}

\setcounter{equation}{0}

In this section, we introduce some random cumulant semigroups, which generalize that of a classical CB-process. Those semigroups are important tools in the study of the CBRE-process. Let $I\subset \mbb{R}$ be an interval and $\zeta = \{\zeta(t): t\in I\}$ a c\`adl\`ag function. Let $\phi$ be a \textit{branching mechanism} given by
 \beqlb\label{eq2.1}
\phi(z) = bz + cz^2 + \int_0^\infty (e^{-uz}-1+uz) m(du), \qquad z\geq 0.
 \eeqlb
where $b\in \mbb{R}$ and $c\ge 0$ are constants and $(z\land z^2)m(dz)$ is a finite measure on $(0,\infty)$. From Theorem~6.10 in Li (2011) it follows that, for any $t\in I$ and $\lambda\ge 0$, there is a unique positive solution $r\mapsto u^\zeta_{r,t}(\lambda)$ to the integral evolution equation
 \beqlb\label{eq2.2}
u^\zeta_{r,t}(\lambda) = \lambda - \int_r^t e^{\zeta(s)}\phi(e^{-\zeta(s)}u^\zeta_{s,t}(\lambda))ds, \qquad r\in I\cap (-\infty, t].
 \eeqlb
Moreover, there is an inhomogeneous transition semigroup $(P^\zeta_{r,t})_{t\ge r\in I}$ on $[0,\infty)$ defined by
 \beqlb\label{eq2.3}
\int_{[0,\infty)} e^{-\lambda y} P^\zeta_{r,t}(x,dy)
 =
e^{-xu^\zeta_{r,t}(\lambda)}, \qquad \lambda\ge 0.
 \eeqlb
By a simple transformation, we can define another inhomogeneous transition semigroup $(Q^\zeta_{r,t})_{t\ge r\in I}$ on $[0,\infty)$ by
 \beqlb\label{eq2.4}
\int_{[0,\infty)} e^{-\lambda y} Q^\zeta_{r,t}(x,dy)
 =
e^{-xv^\zeta_{r,t}(\lambda)}, \qquad \lambda\ge 0,
 \eeqlb
where
 \beqlb\label{eq2.5}
v^\zeta_{r,t}(\lambda) = e^{-\zeta(r)} u^\zeta_{r,t}(e^{\zeta(t)}\lambda).
 \eeqlb
The uniqueness of the solution to (\ref{eq2.2}) implies that
 \beqlb\label{eq2.6}
u^\zeta_{r,t}(\lambda) = u^\zeta_{r,s}\circ u^\zeta_{s,t}(\lambda), \qquad \lambda\ge 0, t\ge s\ge r\in I.
 \eeqlb
There is a similar relation for $(v^\zeta_{r,t})_{t\ge r\in I}$. By Lebesgue's theorem one can see $r\mapsto u^\zeta_{r,t}(\lambda)$ is also the unique positive continuous solution to the differential equation
 \beqlb\label{eq2.7}
\frac{d}{dr}u^\zeta_{r,t}(\lambda) = e^{\zeta(r)}\phi(e^{-\zeta(r)}u^\zeta_{r,t}(\lambda)), \qquad \mbox{a.e.~} r\in I\cap (-\infty,t]
 \eeqlb
with terminal condition $u^\zeta_{t,t}(\lambda) = \lambda$.

\begin{proposition}\label{th2.1} Let $(P^\zeta_{r,t})_{t\ge r\in I}$ be defined by (\ref{eq2.3}). Then for any $x\ge 0$ and $t\ge r\in I$ we have
 \beqnn
\int_{[0,\infty)} y P^\zeta_{r,t}(x,dy) = xe^{-b(t-r)}.
 \eeqnn
\end{proposition}

\proof By differentiating both sides of (\ref{eq2.2}) and solving the resulted integral equation we obtain $(d/d\lambda)u^\zeta_{r,t}(0+) = e^{-b(t-r)}$. Then we get the desired equality by differentiating both sides of (\ref{eq2.3}). \qed

\begin{proposition}\label{th2.2} If $b\ge 0$, then $t\mapsto u^\zeta_{r,t}(\lambda)$ is decreasing on $I\cap [r,\infty)$ and $r\mapsto u^\zeta_{r,t}(\lambda)$ is increasing on $I\cap (-\infty, t]$. \end{proposition}

\proof From (\ref{eq2.3}) we see that $\lambda\mapsto u^\zeta_{r,t}(\lambda)$ is increasing. Since $b\ge 0$, we have $\phi(z)\ge 0$ for every $z\ge 0$. Then (\ref{eq2.2}) implies $u^\zeta_{r,t}(\lambda)\le \lambda$. By (\ref{eq2.6}) we see $u^\zeta_{r,t}(\lambda)\le u^\zeta_{s,t}(\lambda)$ and $u^\zeta_{r,t}(\lambda)\le u^\zeta_{r,s}(\lambda)$ for $r\le s\le t\in I$. \qed

When the function $\zeta$ is degenerate ($\zeta(t) = 0$ for all $t\in I$), both $(u^\zeta_{r,t})_{t\ge r\in I}$ and $(v^\zeta_{r,t})_{t\ge r\in I}$ reduce to the cumulant semigroup of a classical CB-process with branching mechanism $\phi$; see, e.g., Chapter~3 of Li (2011). In the general case, we may think of $(u^\zeta_{r,t})_{t\ge r\in I}$ as an \textit{inhomogeneous cumulant semigroup} determined by the \textit{time-dependent branching mechanism} $(s,z)\mapsto e^{\zeta(s)} \phi(e^{-\zeta(s)}z)$. The idea of the proof of Theorem~6.10 in Li (2011) is to reduce the construction of an inhomogeneous cumulant semigroup to that of a homogeneous one by some time-space processes. The transformation from $(u^\zeta_{r,t})_{t\ge r\in I}$ to $(v^\zeta_{r,t})_{t\ge r\in I}$ is a time-dependent variation of the one used in the proof of Theorem~6.1 in Li (2011).

We next consider some randomization of the inhomogeneous cumulant semigroups defined above. Let $(\Omega, \mcr{F},\mcr{F}_t,\mbf{P})$ be a filtered probability space satisfying the usual hypotheses. Let $a\in \mbb{R}$ and $\sigma\geq 0$ be given constants and $(1\wedge z^2)\nu(dz)$ a finite measure on $(0,\infty)$. Suppose that $\{B(t): t\ge 0\}$ is an $(\mcr{F}_t)$-Brownian motion and $N(ds,dz)$ is an $(\mcr{F}_t)$-Poisson random measure on $(0,\infty)\times \mbb{R}$ with intensity $ds\nu(dz)$. Let $\{\xi(t): t\ge 0\}$ be a $(\mcr{F}_t)$-L\'evy process with the following L\'evy-It\^o decomposition:
 \beqlb\label{eq2.8}
\xi(t) = \xi(0) + at + \sigma B(t) + \int_0^t \int_{[-1,1]} z \tilde{N}(ds,dz) + \int_0^t \int_{[-1,1]^c} z N(ds,dz),
 \eeqlb
where $[-1,1]^c = \mbb{R}\setminus [-1,1]$. Let $u^\xi_{r,t}(\lambda)$ and $v^\xi_{r,t}(\lambda)$ be defined by (\ref{eq2.2}) and (\ref{eq2.5}) with $\zeta=\xi$. From (\ref{eq2.2}) we see that $r\mapsto v^\xi_{r,t}(\lambda)$ is the unique positive solution to
 \beqlb\label{eq2.9}
v^\xi_{r,t}(\lambda) = e^{\xi(t)-\xi(r)}\lambda - \int_r^t e^{\xi(s)-\xi(r)}\phi(v^\xi_{s,t}(\lambda)) ds, \qquad 0\le r\le t.
 \eeqlb
Let $\{L(t): t\ge 0\}$ be the $(\mcr{F}_t)$-L\'evy process with L\'evy-It\^o decomposition:
 \beqlb\label{eq2.10}
L(t) \ar=\ar L(0) + \beta t + \sigma B(t) + \int_0^t \int_{[-1,1]} (e^z-1)\tilde{N}(ds,dz) \cr
  \ar\ar\qquad\qquad\qquad
+ \int_0^t \int_{[-1,1]^c} (e^z-1) N(ds,dz),
 \eeqlb
where $L(0)=\xi(0)$ and
 \beqlb\label{eq2.11}
\beta = a + \frac{\sigma^2}{2} + \int_{[-1,1]} (e^z-1-z) \nu(dz).
 \eeqlb
Clearly, the two processes $\{\xi(t): t\ge 0\}$ and $\{L(t): t\ge 0\}$ generate the same filtration. From (\ref{eq2.9}) we see that the left-continuous process $\{v^\xi_{t-s,t}(\lambda): 0\le s\le t\}$ is progressively measurable with respect to the filtration generated by the L\'evy process $\{L_t(s) := L(t-)-L((t-s)-): 0\le s\le t\}$.

\begin{theorem}\label{th2.3} For any $t\ge 0$ and $\lambda\ge 0$, the process $\{v^\xi_{r,t}(\lambda): 0\le r\le t\}$ is the pathwise unique positive solution to
 \beqlb\label{eq2.12}
v^\xi_{r,t}(\lambda) = \lambda - \int_r^t \phi(v^\xi_{s,t}(\lambda))ds + \int_r^t v^\xi_{s,t}(\lambda)L(\overleftarrow{ds}), \qquad 0\le r\le t,
 \eeqlb
where the backward stochastic integral is defined by
 \beqnn
\int_r^t v^\xi_{s,t}(\lambda)L(\overleftarrow{ds})
 =
\int_0^{(t-r)-} v^\xi_{t-s,t}(\lambda)L_t(ds).
 \eeqnn
\end{theorem}

\proof Let $\xi_t(r) = \xi(t-) - \xi((t-r)-)$ and $B_t(r) = B(t) - B(t-r)$ for $0\le r\le t$. Let $N_t(ds,dz)$ be the Poisson random measure defined by
 \beqnn
N_t([0,r]\times B) = N([t-r,t]\times B), \qquad 0\le r\le t, B\in \mcr{B}(\mbb{R}).
 \eeqnn
From (\ref{eq2.8}) we have
 \beqnn
\xi_t(r) = ar + \sigma B_t(r) + \int_0^r\int_{[-1,1]} z \tilde{N}_t(ds,dz) + \int_0^r \int_{[-1,1]^c} z N_t(ds,dz).
 \eeqnn
On the other hand, from (\ref{eq2.9}) we have $f_t(r) := v^\xi_{(t-r)-,t}(\lambda) = e^{\xi_t(r)} F_t(r)$, where
 \beqnn
F_t(r) = \lambda - \int_0^r e^{-\xi_t(s)}\phi(f_t(s)) ds.
 \eeqnn
By It\^o's formula,
 \beqnn
f_t(r) \ar=\ar \lambda + \int_0^r e^{\xi_t(s-)}F_t(s) \xi_t(ds) + \frac{\sigma^2}{2}\int_0^r e^{\xi_t(s)} F_t(s) ds \cr
 \ar\ar
+ \int_0^r\int_{\mbb{R}} e^{\xi_t(s-)} F_t(s)(e^z-1-z) N_t(ds,dz) + \int_0^r e^{\xi_t(s)} F_t(ds) \cr
 \ar=\ar
\lambda + \int_0^r f_t(s-)L_t(ds) - \int_0^r \phi(f_t(s))ds.
 \eeqnn
It follows that
 \beqnn
v^\xi_{r,t}(\lambda) \ar=\ar \lambda + \int_0^{r-} f_t(s-)L_t(ds) - \int_0^r \phi(f_t(s))ds \cr
 \ar=\ar
\lambda + \int_{0}^{(t-r)-} v^\xi_{t-s,t}(\lambda)L_t(ds) - \int_0^{t-r} \phi(v^\xi_{t-s,t}(\lambda))ds \cr
 \ar=\ar
\lambda + \int_r^t v^\xi_{s,t}(\lambda)L(\overleftarrow{ds}) - \int_r^t \phi(v^\xi_{t-s,t}(\lambda))ds.
 \eeqnn
That proves the existence of the solution to (\ref{eq2.12}). Conversely, assuming $r\mapsto v^\xi_{r,t}(\lambda)$ is a solution to (\ref{eq2.12}), one can use similar calculations to see it also solves (\ref{eq2.9}). Then the pathwise uniqueness for (\ref{eq2.12}) is a consequence of that for (\ref{eq2.9}). \qed


\section{Construction of CBRE-processes}

\setcounter{equation}{0}

Let $(\Omega, \mcr{F},\mcr{F}_t,\mbf{P})$ be a filtered probability space satisfying the usual hypotheses. Let $\{\xi(t): t\ge 0\}$ and $\{L(t): t\ge 0\}$ be $(\mcr{F}_t)$-L\'evy processes given as in Section~2. Let $b\in \mbb{R}$ and $c\geq 0$ be constants and $(z\wedge z^2)m(dz)$ a finite measure on $(0,\infty)$. Suppose that $\{W(t): t\ge 0\}$ is another $(\mcr{F}_t)$-Brownian motion and $M(ds,dz,du)$ is an $(\mcr{F}_t)$-Poisson random measure on $(0,\infty)^3$ with intensity $dsm(dz)du$. We assume both of those are independent of the L\'evy process $\{L(t): t\ge 0\}$. Given a positive $\mcr{F}_0$-measurable random variable $X(0)$, we consider the following stochastic integral equation:
 \beqlb\label{eq3.1}
X(t) \ar=\ar X(0) - b\int_0^t X(s)ds + \int_0^t\sqrt{2cX(s)}dW(s) \cr
 \ar\ar
+ \int_0^t \int_0^\infty\int_0^{X(s-)} z \tilde{M}(ds,dz,du) + \int_0^t X(s-)dL(s),
 \eeqlb
where $\tilde{M}(ds,dz,du) = M(ds,dz,du) - dsm(dz)du$.

\begin{theorem}\label{th3.1} There is a unique positive strong solution $\{X(t): t\ge 0\}$ to (\ref{eq3.1}).
\end{theorem}

\proof Let $E = \{1,2\}$ and $U_0 = A_0\cup B_0$, where $A_0 = \{1\}\times (0,\infty)^2$ and $B_0 = \{2\}\times (0,\infty)$. Let $\pi(du) = \delta_1(du) + \delta_2(du)$ for $u\in E$. Then $W(ds,du) := dW(s)\delta_1(du) + dB(s)\delta_2(du)$ is a Gaussian white noise on $(0,\infty)\times E$ with intensity $ds\pi(dz)$. Let $\mu_0(dy,dz,du) = \delta_1(dy)m(dz)du$ for $(y,z,u)\in A_0$ and $\mu_0(dy,dz) = \delta_2(dy)\nu(dz)$ for $(y,z)\in B_0$. Then $N_0(ds,dy,dz,du) := \delta_1(dy)M(ds,dz,du)$ is a Poisson random measure on $(0,\infty)\times A_0$ with intensity $ds\delta_1(dy)m(dz)du$ and $N_0(ds,dy,dz) := \delta_2(dy)N(ds,dz)$ is a Poisson random measure on $(0,\infty)\times B_0$ with intensity $ds\delta_2(dy)\nu(dz)$. Let $b(x) = (\beta-b)x$ for $x\in [0,\infty)$ and $\sigma(x,u) = \sqrt{2cx}1_{\{u=1\}} + \sigma x1_{\{u=2\}}$ for $(x,u)\in [0,\infty)\times E$. Let $g_0(x,y,z,u) = z1_{\{u\le x\}}$ for $(x,y,z,u)\in [0,\infty)\times A_0$ and $g_0(x,y,z) = x(e^z-1)$ for $(x,y,z)\in [0,\infty)\times B_0$. By Theorem~2.5 in Dawson and Li (2012), there is a unique positive strong solution to the stochastic equation
 \beqnn
X(t) \ar=\ar X(0) + \int_0^t b(X(s))ds + \int_0^t\int_{A_0} g_0(X(s-),y,z,u)\tilde{N}_0(ds,dy,dz,du) \cr
 \ar\ar
+ \int_0^t\int_{E} \sigma(X(s),u)W(ds,du) + \int_0^t\int_{B_0} g_0(X(s-),y,u)\tilde{N}_0(ds,dy,du).
 \eeqnn
The above equation can be rewritten into
 \beqnn
X(t) \ar=\ar X(0) + \int_0^t\sqrt{2cX(s)}dW(s) + \int_0^t \int_0^\infty\int_0^{X(s-)} z \tilde{M}(ds,dz,du) \cr
 \ar\ar
+\, (\beta-b)\int_0^t X(s)ds + \sigma\int_0^t X(s)dB(s) + \int_0^t \int_{[-1,1]} X(s-)(e^z-1)\tilde{N}(ds,dz).
 \eeqnn
Since the process
 \beqnn
t\mapsto \int_0^t \int_{[-1,1]^c} (e^z-1) N(ds,dz)
 \eeqnn
has at most a finite number of jumps in each bounded time interval, as in the proof of Proposition~2.2 in Fu and Li (2010), one can see that there is also a pathwise unique positive strong solution to
 \beqlb\label{eq3.2}
X(t) \ar=\ar X(0) + \int_0^t\sqrt{2cX(s)}dW(s) + \int_0^t \int_0^\infty\int_0^{X(s-)} z \tilde{M}(ds,dz,du) \cr
 \ar\ar
+\, (\beta-b)\int_0^t X(s)ds + \sigma\int_0^t X(s)dB(s) + \int_0^t \int_{[-1,1]} X(s-)(e^z-1) \tilde{N}(ds,dz) \cr
 \ar\ar
+ \int_0^t \int_{[-1,1]^c} X(s-)(e^z-1) N(ds,dz).
 \eeqlb
The above equation is just a reformulation of (\ref{eq3.1}). Then we have the result of the theorem. \qed

\begin{remark}\label{th3.2} Since the stochastic equations in Section~2 of Dawson and Li (2012) are formulated for Gaussian and Poisson noises in abstract space, they are quite flexible for applications. For example, the proof of Theorem~\ref{th3.1} given above can be modified to more complex models with extra structures such as immigration, competition and so on. \end{remark}

We call the solution $\{X(t): t\ge 0\}$ to (\ref{eq3.1}) a \textit{CBRE-process}, which is a c\`adl\`ag strong Markov process. Here the random environment is provided by the L\'evy process $\{L(t): t\ge 0\}$. By It\^o's formula one can see $\{X(t): t\ge 0\}$ has strong generator $A$ defined as follows: For $f\in C^2(\mbb{R}_+)$,
 \beqlb\label{eq3.3}
Af(x) \ar=\ar (\beta-b) xf'(x) + cxf''(x) + x\int_0^\infty[f(x+z)-f(x)-zf'(x)]m(dz)\cr
 \ar\ar
+ \frac{\sigma^2}{2}x^2f''(x) + \int_{[-1,1]} [f(xe^z) - f(x) - x(e^z-1)f'(x)] \nu(dz)\cr
 \ar\ar
+ \int_{[-1,1]^c} [f(xe^z)-f(x)] \nu(dz).
 \eeqlb

\begin{proposition}\label{th3.3} Let $Z(t) = X(t)e^{-\xi(t)}$ for $t\ge 0$. Then we have
 \beqlb\label{eq3.4}
Z(t) \ar=\ar X(0) - b\int_0^t e^{-\xi(s)}X(s) ds + \int_0^t e^{-\xi(s)}\sqrt{2cX(s)} dW(s) \cr
 \ar\ar\qquad
+ \int_0^t\int_0^\infty\int_0^{X(s-)} ze^{-\xi(s-)} \tilde{M}(ds,dz,du).
 \eeqlb
In particular, the process $\{Z(t): t\ge 0\}$ is a positive local martingale when $b=0$.
\end{proposition}

\proof Let $f(x,y) = xe^{-y}$. Then $xf'_x(x,y) = -f'_y(x,y) = f(x,y)$ and $-xf''_{xy}(x,y) = f''_{yy}(x,y) = f(x,y)$. Observe that the Poison random measure $N(ds,dz)$ actually does not produce any jump of $t\mapsto Z(t)$. By (\ref{eq3.1}) and It\^{o}'s formula,
 \beqnn
Z(t) \ar=\ar X(0) - b\int_0^t f'_x(X(s),\xi(s))X(s)ds + \int_0^t f'_x(X(s),\xi(s))\sqrt{2cX(s)} dW(s) \cr
 \ar\ar
+ \int_0^t\int_0^\infty\int_0^{X(s-)} f'_x(X(s-),\xi(s-))z \tilde{M}(ds,dz,du) \cr
 \ar\ar
+\, \beta\int_0^t f'_x(X(s),\xi(s))X(s) ds + \sigma\int_0^t f'_x(X(s),\xi(s))X(s) dB((s) \cr
 \ar\ar
+ \int_0^t\int_{[-1,1]} f'_x(X(s-),\xi(s-))X(s-)(e^z-1) \tilde{N}(ds,dz) \cr
 \ar\ar
+ \int_0^t\int_{[-1,1]^c} f'_x(X(s-),\xi(s-))X(s-)(e^z-1) N(ds,dz) \cr
 \ar\ar
+\, a\int_0^t f'_y(X(s),\xi(s))ds + \int_0^t\int_{[-1,1]} f'_y(X(s-),\xi(s-))z \tilde{N}(ds,dz) \cr
 \ar\ar
+\, \sigma\int_0^t f'_y(X(s),\xi(s)) dB((s)
+ \int_0^t\int_{[-1,1]^c} f'_y(X(s-),\xi(s-))z N(ds,dz) \cr
 \ar\ar
+\, \sigma^2\int_0^t \Big[f''_{xy}(X(s),\xi(s))X(s) + \frac{1}{2}f''_{yy}(X(s),\xi(s))\Big] ds \cr
 \ar\ar
+ \int_0^t\int_0^\infty\int_0^{X(s-)} [f(X(s-)+z,\xi(s)) - f(X(s-),\xi(s)) \cr
 \ar\ar\hskip6cm
- f'_x(X(s-),\xi(s))z] M(ds,dz,du) \cr
\ar\ar
+ \int_0^t\int_{\mbb{R}} [f(X(s-)e^z,\xi(s-)+z) - f(X(s-),\xi(s-)) \cr
 \ar\ar\hskip2cm
- f'_x(X(s-),\xi(s-))X(s-)(e^z-1) - f'_y(X(s-),\xi(s-))z] N(ds,dz) \cr
 \ar=\ar
X(0) - b\int_0^t f'_x(X(s),\xi(s))X(s)ds + \int_0^t f'_x(X(s),\xi(s))\sqrt{2cX(s)} dW(s) \cr
 \ar\ar
+ \int_0^t\int_0^\infty\int_0^{X(s-)} f'_x(X(s-),\xi(s-))z \tilde{M}(ds,dz,du) \cr
 \ar\ar
+\, \beta\int_0^t f'_x(X(s),\xi(s))X(s) ds + \sigma\int_0^t f'_x(X(s),\xi(s))X(s) dB((s) \cr
 \ar\ar
+\, a\int_0^t f'_y(X(s),\xi(s))ds + \sigma\int_0^t f'_y(X(s),\xi(s)) dB((s) \cr
 \ar\ar
+\, \sigma^2\int_0^t \Big[f''_{xy}(X(s),\xi(s))X(s) + \frac{1}{2}f''_{yy}(X(s),\xi(s))\Big] ds \cr
 \ar\ar
- \int_0^tds \int_{[-1,1]}[f'_x(X(s),\xi(s))X(s)(e^z-1) + f'_y(X(s),\xi(s))z]\nu(dz).
 \eeqnn
By reorganizing the terms on the right-hand side we get the desired equality. \qed

It is easy to see that the two L\'evy processes $\{\xi(t): t\ge 0\}$ and $\{L(t): t\ge 0\}$ generate the same $\sigma$-algebra. Let $\mbf{P}^\xi$ denote the quenched law given $\{\xi(t): t\ge 0\}$ or $\{L(t): t\ge 0\}$.

\begin{theorem}\label{th3.4} Let $(P^\xi_{r,t})_{t\ge r\ge 0}$ and $(Q^\xi_{r,t})_{t\ge r\ge 0}$ be defined by (\ref{eq2.3}) and (\ref{eq2.4}), respectively, with $\zeta=\xi$. Then for any $\lambda\geq 0$ and $t\geq r\geq 0$ we have
 \beqlb\label{eq3.5}
\mbf{P}^\xi [e^{-\lambda Z(t)}|\mcr{F}_r] = \exp\{-Z(r)u^\xi_{r,t}(\lambda)\}
 =
\int_{[0,\infty)} e^{-\lambda y} P_{r,t}^\xi(Z(r),dy)
 \eeqlb
and
 \beqlb\label{eq3.6}
\mbf{P}^\xi [e^{-\lambda X(t)}|\mcr{F}_r]
 =
\exp\{-X(r)v^\xi_{r,t}(\lambda)\}
 =
\int_{[0,\infty)} e^{-\lambda y} Q_{r,t}^\xi(X(r),dy).
 \eeqlb
\end{theorem}

\proof Fix $\lambda\geq 0$. For $t\ge r\ge 0$, let $H_t(r) = \exp\{-Z(r)u^\xi_{r,t}(\lambda)\}$. By (\ref{eq2.2}) and Proposition~\ref{th3.3}, given the environment $\{\xi(t): t\ge 0\}$, we can use It\^{o} formula to see
 \beqnn
H_t(t)
 \ar=\ar
H_t(r) - \int_r^t H_t(s)Z(s)e^{\xi(s)}\phi(e^{-\xi(s)}u^\xi_{s,t}(\lambda))ds \cr
 \ar\ar
-\, \int_r^t H_t(s-)u^\xi_{s,t}(\lambda)dZ(s) + c\int_r^t H_t(s)u^\xi_{s,t}(\lambda)^2 e^{-2\xi(s)}X(s)ds \cr
 \ar\ar
+ \int_r^t\int_0^\infty \int_0^{X(s-)} H_t(s-) \Big[\exp\{-ze^{-\xi(s-)}u^\xi_{s,t}(\lambda)\} \cr
 \ar\ar\hskip5cm
- 1 + ze^{-\xi(s-)}u^\xi_{s,t}(\lambda)\Big] M(ds,dz,du)\cr
 \ar=\ar
H_t(r) - \int_r^t H_t(s)X(s)\phi(e^{-\xi(s)}u^\xi_{s,t}(\lambda))ds + b\int_r^t H_t(s)u^\xi_{s,t}(\lambda)e^{-\xi(s)}X(s) ds \cr
 \ar\ar
-\, \int_r^t H_t(s)u^\xi_{s,t}(\lambda)e^{-\xi(s)}\sqrt{2cX(s)} dW(s) + c\int_r^t H_t(s)e^{-2\xi(s)} u^\xi_{s,t}(\lambda)^2 X(s)ds \cr
 \ar\ar
+ \int_r^t\int_0^\infty \int_0^{X(s-)} H_t(s-) \Big[\exp\{-ze^{-\xi(s-)}u^\xi_{s,t}(\lambda)\} \cr
 \ar\ar\hskip5cm
- 1 + ze^{-\xi(s-)}u^\xi_{s,t}(\lambda)\Big] M(ds,dz,du) \cr
 \ar=\ar
H_t(r) - \int_r^t H_t(s-)u^\xi_{s,t}(\lambda)e^{-\xi(s)}\sqrt{2cX(s)} dW(s) \cr
 \ar\ar
+ \int_r^t\int_0^\infty \int_0^{X(s-)} H_t(s-) \Big[\exp\{-ze^{-\xi(s-)}u^\xi_{s,t}(\lambda)\} \cr
 \ar\ar\hskip5cm
- 1 + ze^{-\xi(s-)}u^\xi_{s,t}(\lambda)\Big] \tilde{M}(ds,dz,du).
 \eeqnn
Since $\{H_t(r): t\ge r\}$ is a bounded process, by taking the conditional expectation in both sides we get $\mbf{P}^\xi[H_t(t)|\mcr{F}_r] = H_t(r)$. That gives (\ref{eq3.5}), and as a consequence we get (\ref{eq3.6}). \qed

\begin{corollary}\label{th3.5} If $\mbf{P}[Z(0)] = \mbf{P}[X(0)]< \infty$, then $\{e^{bt}Z(t): t\ge 0\}$ is a martingale. \end{corollary}

\proof Let $t\ge r\ge 0$ and let $F$ be a bounded random variable measurable with respect to the $\sigma$-algebra generated by $\mcr{F}_r\cup \sigma(\xi)$. By (\ref{eq3.5}) and Proposition~\ref{th2.1}, we have
 \beqnn
\mbf{P}[Fe^{bt}Z(t)] = \mbf{P}\{Fe^{bt}\mbf{P}^\xi[Z(t)|\mcr{F}_r]\}
 =
\mbf{P}\bigg[Fe^{bt}\int_{[0,\infty)} y P_{r,t}^\xi(Z(r),dy)\bigg]
 =
\mbf{P}[Fe^{br}Z(r)].
 \eeqnn
Then $\{e^{bt}Z(t): t\ge 0\}$ is a martingale. \qed

By Theorem~\ref{th3.4} we see that $\{Z(t): t\ge 0\}$ and $\{X(t): t\ge 0\}$ are actually CB-processes under the quenched law with inhomogeneous cumulant semigroups $(u^\xi_{r,t})_{t\ge r\ge 0}$ and $(v^\xi_{r,t})_{t\ge r\ge 0}$, respectively. The next theorem gives a characterization of the transition semigroup of $\{X(t): t\ge 0\}$ under the annealed law.

\begin{theorem}\label{th3.6} The Markov process $\{X(t): t\ge 0\}$ has Feller transition semigroup $(Q_t)_{t\ge 0}$ defined by
 \beqlb\label{eq3.7}
\int_{[0,\infty)} e^{-\lambda y} Q_t(x,dy)
 =
\mbf{P}[e^{-xv^\xi_{0,t}(\lambda)}], \qquad \lambda\ge 0.
 \eeqlb
\end{theorem}

\proof By (\ref{eq2.4}) one can see that (\ref{eq3.7}) defines a probability kernel $Q_t(x,dy)$. In view of (\ref{eq3.6}), for any bounded $\mcr{F}_r$-measurable random variable $F$ we have
 \beqnn
\mbf{P}[Fe^{-\lambda X(t)}]
 =
\mbf{P}[F\mbf{P}^\xi(e^{-\lambda X(t)}|\mcr{F}_r)]
 =
\mbf{P}[F\exp\{-X(r)v^\xi_{r,t}(\lambda)\}].
 \eeqnn
The pathwise uniqueness of the solution to (\ref{eq2.12}) implies that the random variable $v^\xi_{r,t}(\lambda)$ is measurable with respect to the $\sigma$-algebra generated by $\{L(s)-L(t): r\le s\le t\}$ and is identically distributed with $v^\xi_{0,t-r}(\lambda)$. It follows that
 \beqnn
\mbf{P}[Fe^{-\lambda X(t)}]
 =
\mbf{P}\bigg[F\int_{[0,\infty)} e^{-y\lambda} Q_{t-r}(X(r),dy)\bigg].
 \eeqnn
Then $\{X(t): t\ge 0\}$ has transition semigroup $(Q_t)_{t\ge 0}$. The Feller property is immediate by (\ref{eq3.7}). \qed

Under the annealed law, the process $\{Z(t): t\ge 0\}$ usually does not satisfy the Markov property, but $\{(\xi(t),Z(t)): t\ge 0\}$ is a two-dimensional Markov process. Let $\mcr{F}_\infty = \sigma(\cup_{t\ge 0}\mcr{F}_t)$. By Corollary~\ref{th3.5}, there is a probability measure $\tilde{\mbf{P}}$ on $(\Omega, \mcr{F}_\infty)$ so that $\tilde{\mbf{P}}(F) = \mbf{P}[Fe^{bt}Z(t)]$ for each bounded $\mcr{F}_t$-measurable random variable $F$. Let $\tilde{\mbf{P}}^\xi$ denote the conditional law under $\tilde{\mbf{P}}$ given the environment $\{\xi(t): t\ge 0\}$.

\begin{theorem}\label{th3.7} For any $\lambda\geq 0$ and $t\geq r\geq 0$, we have
 \beqlb\label{eq3.8}
\tilde{\mbf{P}}^\xi[e^{-\lambda X(t)}|\mcr{F}_r]
 =
\exp\bigg\{-X(r)v^\xi_{r,t}(\lambda) - \int_r^t\phi_0'(v^\xi_{s,t}(\lambda)) ds\bigg\},
 \eeqlb
where $\phi_0'(z) = \phi'(z)-b$.
\end{theorem}

\proof From (\ref{eq2.4}) one can see that $\lambda\mapsto v^\xi_{r,t}(\lambda)$ is infinitely differentiable in $(0,\infty)$. By differentiating both sides of (\ref{eq2.2}) we obtain
 \beqnn
\frac{d}{d\lambda}u^\xi_{r,t}(\lambda) = 1 - \int_r^t \phi'(e^{-\xi(s)}u^\xi_{s,t}(\lambda)) \frac{d}{d\lambda}u^\xi_{s,t}(\lambda)ds.
 \eeqnn
Then we can solve the equation to get
 \beqnn
\frac{d}{d\lambda}u^\xi_{r,t}(\lambda)
 =
\exp\bigg\{-\int_r^t \phi'(e^{-\xi(s)}u^\xi_{s,t}(\lambda)) ds\bigg\}.
 \eeqnn
Let $F$ be a bounded random variable measurable with respect to the $\sigma$-algebra generated by $\mcr{F}_r\cup \sigma(\xi)$. From (\ref{eq3.5}) it follows that
 \beqnn
\mbf{P}[Fe^{-\lambda Z(t)}]
 =
\mbf{P}[F\mbf{P}^\xi(e^{-\lambda Z(t)}|\mcr{F}_r)]
 =
\mbf{P}[F\exp\{-Z(r)u^\xi_{r,t}(\lambda)\}].
 \eeqnn
By differentiating both sides in $\lambda>0$ we have
 \beqnn
\mbf{P}[Fe^{-\lambda Z(t)}Z(t)]
 =
\mbf{P}\bigg[F\exp\bigg\{-Z(r)u^\xi_{r,t}(\lambda) - \int_r^t\phi'(e^{-\xi(s)}u^\xi_{s,t}(\lambda)) ds\bigg\}Z(r)\bigg],
 \eeqnn
and hence
 \beqnn
\mbf{P}[Fe^{-\lambda X(t)}Z(t)]
 =
\mbf{P}\bigg[F\exp\bigg\{-X(r)v^\xi_{r,t}(\lambda) - \int_r^t\phi'(v^\xi_{s,t}(\lambda)) ds\bigg\} Z(r)\bigg].
 \eeqnn
It follows that
 \beqnn
\tilde{\mbf{P}}[Fe^{-\lambda X(t)}]
 =
\tilde{\mbf{P}}\bigg[F\exp\bigg\{-X(r)v^\xi_{r,t}(\lambda) - \int_r^t\phi_0'(v^\xi_{s,t}(\lambda)) ds\bigg\}\bigg].
 \eeqnn
Then we have (\ref{eq3.8}). The extension of the equality to $\lambda\ge 0$ is immediate. \qed

By Theorem~\ref{th3.7} one can show as in the proof of Theorem~\ref{th3.6} that $\{X(t): t\ge 0\}$ is a Markov process under $\tilde{\mbf{P}}$ with Feller transition semigroup $(\tilde{Q}_t)_{t\ge 0}$ defined by
 \beqlb\label{eq3.9}
\int_{[0,\infty)} e^{-\lambda y} \tilde{Q}_t(x,dy)
 =
\mbf{P}\bigg[\exp\bigg\{-xv^\xi_{0,t}(\lambda) - \int_0^t\phi_0'(v^\xi_{s,t}(\lambda)) ds\bigg\}\bigg], \qquad \lambda\ge 0.
 \eeqlb
This is a special case of a larger class of transition semigroups to be given in Section~5.


\section{Survival and extinction probabilities}

\setcounter{equation}{0}

In this section, we assume $X(0)=x> 0$ is a deterministic constant for simplicity. Let $\mbf{P}$ or $\mbf{P}_x$ denote the annealed law and $\mbf{P}^\xi$ or $\mbf{P}^\xi_x$ the quenched law given $\{\xi(t): t\ge 0\}$. Let $\tau_0 = \inf\{t\ge 0: X(t)=Z(t)=0\}$ denote the \textit{extinction time} of the CBRE-process. From (\ref{eq2.4}) one can see that $v^\xi_{0,t}(\lambda)$ is increasing in $\lambda\ge 0$. For $t>0$ let $\bar{v}^\xi_{0,t} := \lim_{\lambda\to \infty} v^\xi_{0,t}(\lambda) \in[0,\infty]$. Then
 \beqlb\label{eq5.1}
\bar{u}^\xi_{0,t} := \lim_{\lambda\to \infty}u^\xi_{0,t}(\lambda)
 =
\lim_{\lambda\to \infty}v^\xi_{0,t}(e^{-\xi(t)}\lambda) = \bar{v}^\xi_{0,t}.
 \eeqlb
By (\ref{eq3.5}) and (\ref{eq3.6}) we have the following characterizations of the extinction probabilities:
 \beqlb\label{eq5.2}
\mbf{P}^\xi_x(\tau_0\le t) = \mbf{P}^\xi_x(Z(t)=0) = \mbf{P}^\xi_x(X(t)=0)
 =
e^{-x\bar{u}^\xi_{0,t}} = e^{-x\bar{v}^\xi_{0,t}}
 \eeqlb
and
 \beqlb\label{eq5.3}
\mbf{P}_x(\tau_0\le t) = \mbf{P}_x(Z(t)=0) = \mbf{P}_x(X(t)=0)
 =
\mbf{P}(e^{-x\bar{u}^\xi_{0,t}}) = \mbf{P}(e^{-x\bar{v}^\xi_{0,t}}).
 \eeqlb
We say the branching mechanism $\phi$ satisfies \textit{Grey's condition} if
 \beqlb\label{eq5.4}
\int_1^\infty \phi(z)^{-1}dz< \infty.
 \eeqlb

\begin{theorem}\label{th5.1} The following statements are equivalent:
 \benumerate

\itm[(1)] $\phi$ satisfies Grey's condition;

\itm[(2)] $\mbf{P}_x({Z(t)=0}) = \mbf{P}_x({X(t)=0})> 0$ for some and hence all $t> 0$;

\itm[(3)] $\mbf{P}(\bar{u}^\xi_{0,t}< \infty) = \mbf{P}(\bar{v}^\xi_{0,t}< \infty)> 0$ for some and hence all $t> 0$;

\itm[(4)] $\mbf{P}(\bar{u}^\xi_{0,t}< \infty) = \mbf{P}(\bar{v}^\xi_{0,t}< \infty) = 1$ for some and hence all $t> 0$.

 \eenumerate
\end{theorem}

\proof From (\ref{eq5.3}) we see (2)$\Leftrightarrow$(3)$\Leftarrow$(4). Then we only need to show (3)$\Rightarrow$(1)$\Rightarrow$(4). From (\ref{eq2.7}) we have
 \beqlb\label{eq5.5}
t = \int_0^t\frac{1}{e^{\xi(s)}\phi(e^{-\xi(s)}u^\xi_{s,t}(\lambda))} du^\xi_{s,t}(\lambda).
 \eeqlb
Suppose that $\mbf{P}(\bar{v}^\xi_{0,t}< \infty) = \mbf{P}(\bar{u}^\xi_{0,t}< \infty)> 0$ for some $t> 0$. Choose the constants $0< M_1< M_2< \infty$ so that the event $A := \{\bar{u}^\xi_{0,t}< \infty\}\cap \{M_1\leq e^{\xi(s)}\leq M_2$ for $s\in [0,t]\}$ has strictly positive probability. Since $z\mapsto \phi(z)$ is increasing, by (\ref{eq5.5}) we have on $A$ that
 \beqnn
t \geq\int_0^t\frac{1}{M_2\phi(M_1^{-1}u^\xi_{s,t}(\lambda))} du^\xi_{s,t}(\lambda)
 =
\frac{M_1}{M_2}\int_{M_1^{-1}u^\xi_{0,t}(\lambda)}^{M_1^{-1}\lambda}\frac{1}{\phi(z)} dz.
 \eeqnn
By letting $\lambda\to \infty$ we have on $A$ that
 \beqnn
\frac{M_1}{M_2}\int_{M_2\bar{u}^\xi_{0,t}}^\infty\frac{1}{\phi(z)} dz\leq t.
 \eeqnn
Then (\ref{eq5.4}) holds. That proves (3)$\Rightarrow$(1). Now suppose that Grey's condition (\ref{eq5.4}) is satisfied. Fix any $t>0$. Choose sufficiently large $n\ge 1$ so that the event $\Omega_n := \{1/n\leq e^{\xi(s)}\leq n$ for $s\in [0,t]\}$ has strictly positive probability. By (\ref{eq5.5}), on the event $\Omega_n$ we have
 \beqnn
t \leq\int_0^t\frac{n}{\phi(n^{-1}u^\xi_{s,t}(\lambda))} du^\xi_{s,t}(\lambda)
 =
n^2\int_{n^{-1}u^\xi_{0,t}(\lambda)}^{n^{-1}\lambda}\frac{1}{\phi(z)} dz,
 \eeqnn
which implies
 \beqnn
n^2\int_{n^{-1}\bar{u}^\xi_{0,t}}^\infty\frac{1}{\phi(z)} dz\geq t.
 \eeqnn
It follows that $\bar{u}^\xi_{0,t} = \bar{v}^\xi_{0,t}< \infty$ on $\Omega_n$. Since $\mbf{P}(\cup_{n=1}^\infty \Omega_n) = 1$, we have $\mbf{P}(\bar{v}^\xi_{0,t}< \infty) = \mbf{P}(\bar{u}^\xi_{0,t}< \infty) = 1$. \qed

\begin{corollary}\label{th5.2} Under Grey's condition, for any $t>0$, the function $r\mapsto u(r) := \bar{u}^\xi_{r,t} = \bar{v}^\xi_{r,t}$ on $[0,t)$ is the minimal positive continuous solution to
 \beqlb\label{eq5.6}
\frac{d}{dr}u(r) = e^{\xi(r)}\phi(e^{-\xi(r)}u(r)), \qquad \mbox{a.e.~} r\in (0,t)
 \eeqlb
with terminal condition $u(t-) = \infty$.
\end{corollary}

\proof For any $t>s>r>0$ we have $\bar{u}^\xi_{r,t} = \lim_{\lambda\to \infty} u^\xi_{r,s}(u^\xi_{s,t}(\lambda)) = u^\xi_{r,s}(\bar{u}^\xi_{s,t})$. From (\ref{eq2.7}) we see the differential equation in (\ref{eq5.6}) is satisfied first for a.e.\ $r\in (0,s)$ and then for a.e.\ $r\in (0,t)$. Since $\bar{u}^\xi_{t-,t}\ge u^\xi_{t-,t}(\lambda) = \lambda$ for any $\lambda\ge0$, we have the terminal property $\bar{u}^\xi_{t-,t} = \infty$. Now suppose that $r\mapsto w(r)$ is another positive continuous solution to (\ref{eq5.6}). By the uniqueness of the solution to (\ref{eq2.7}) we have $w(r) = u_{r,s}^\xi(w(s))$ for $0\le r\le s<t$. For any $\lambda\ge 0$ we can choose $s\in (r,t)$ so that $w(s)\ge \lambda e^{|b|t}$. By Proposition~\ref{th2.1} and Jensen's inequality one can see $u^\xi_{s,t}(\lambda)\le \lambda e^{-b(t-s)}\le \lambda e^{|b|t}$. From monotonicity of $\lambda\mapsto u^\xi_{r,s}(\lambda)$ we get $w(r) = u_{r,s}^\xi(w(s))\ge u^\xi_{r,s}(\lambda e^{|b|t})\ge u^\xi_{r,s}(u^\xi_{s,t}(\lambda)) = u^\xi_{r,t}(\lambda)$. Then $w(r)\ge \bar{u}^\xi_{r,t} = \lim_{\lambda\to \infty} u^\xi_{r,t}(\lambda)$. \qed

\begin{theorem}\label{th5.3} Let $\bar{v}^\xi := \downarrow\lim_{t\to \infty} \bar{v}^\xi_{0,t}\in[0,\infty]$ and $\tau_0 := \inf\{t\geq 0: X(t)=0\}$. Then
 $$
\mbf{P}_x(\tau_0< \infty)
 =
\lim_{t\to \infty} \mbf{P}_x(\tau_0\le t)
 =
\lim_{t\to \infty} \mbf{P}_x(X(t)=0)
 =
\mbf{P}[e^{-x \bar{v}^\xi}].
 $$
Moreover, we have $\bar{v}^\xi<\infty$ if and only if Grey's condition (\ref{eq5.4}) holds. \end{theorem}

\proof From $\mbf{P}^\xi_x(\tau_0\le t) = \mbf{P}^\xi_x({X(t)=0}) = e^{-x\bar{v}^\xi_{0,t}}$ we see $t\mapsto \bar{v}^\xi_{0,t}$ is decreasing, so $\bar{v}^\xi$ is well defined. From (\ref{eq5.2}) it follows that
 \beqnn
\mbf{P}_x(\tau_0<\infty) = \lim_{t\to \infty} \mbf{P}_x(\tau_0\le t)
 =
\lim_{t\to \infty} \mbf{P}_x(X(t)=0)
 =
\lim_{t\to \infty}\mbf{P}[e^{-x\bar{v}^\xi_{0,t}}] = \mbf{P}[e^{-x \bar{v}^\xi}].
 \eeqnn
The second statement follows immediately from Theorem~\ref{th5.1}. \qed

\begin{corollary}\label{th5.4} Suppose that $\liminf_{t\to \infty} (\xi(t)-bt) = -\infty$ and Grey's condition (\ref{eq5.4}) holds. Then
 \beqnn
\mbf{P}_x(\tau_0=\infty)
 =
\lim_{t\to \infty} \mbf{P}_x(\tau_0> t)
 =
\lim_{t\to \infty} \mbf{P}_x(X(t)>0) = 0.
 \eeqnn
\end{corollary}

\proof Suppose that $p(x) := \mbf{P}_x(\tau_0=\infty) = 1 - \mbf{P}[e^{-x\bar{v}^\xi}]> 0$. Then we have $\mbf{P}(\bar{v}^\xi>0)> 0$, so $x\mapsto p(x) $ is strictly increasing. Under the assumption $\liminf_{t\to \infty} (\xi(t)-bt) = -\infty$, we have $\liminf_{t\to \infty}X(t) = 0$ as observed in Corollary~2 of Bansaye (2013). Then the stopping time $\sigma = \inf\{t>0: \xi(t)<x/2\}$ is a.s. finite. By Theorem~\ref{th5.3} and the strong Markov property, we have
 \beqnn
p(x) = \mbf{P}_x[\mbf{P}_{X(\sigma)}(\tau_0=\infty)]
 =
\mbf{P}_x[\mbf{P}_{x/2}(\tau_0=\infty)]
 =
1-\mbf{P}[e^{-x\bar{v}^\xi/2}] = p(x/2),
 \eeqnn
which yields a contradiction. Then we must have $p(x) = \mbf{P}_x(\tau_0=\infty)=0$. \qed

\begin{theorem}\label{th5.5} Under Grey's condition, the transition semigroup $(Q_t)_{t\geq 0}$ defined by (\ref{eq3.7}) has the strong Feller property. \end{theorem}

\proof We here need a construction of the CBRE-process for all initial values. Let $W(ds,du)$ be a time-space Gaussian white noise on $(0,\infty)^2$ with intensity $dsdu$. By a modification of the proof of Theorem~\ref{th3.1}, one can see for each $x\ge 0$ there is a unique positive strong solution to the stochastic equation
 \beqnn
Y_t(x)\ar=\ar x - b\int_0^t Y_s(x) ds + \sqrt{2c}\int_0^t\int_0^{Y_{s-}(x)}W(ds,du) \cr
 \ar\ar
+ \int_0^t\int_0^{\infty}\int_0^{Y_{s-}(x)}z\tilde{M}(ds,dz,du)+ \int_0^t Y_{s-}(x)L(ds).
 \eeqnn
Clearly, the process $\{Y_t(x): t\ge 0\}$ is equivalent to the solution $\{X(t): t\ge 0\}$ to (\ref{eq3.1}) with $X(0)=x$. As in the proof of Theorem~3.2 in Dawson and Li (2012), one can show that, for $y\ge x\ge 0$ and $t\ge0$ we have $Y_t(y)\ge Y_t(x)$ and $\{Y_t(y)-Y_t(x): t\ge 0\}$ is equivalent to the solution $\{X(t): t\ge 0\}$ to (\ref{eq3.1}) with $X(0)=y-x$. Now let $f\in b\mathscr{B}(\mathbb{R}_+)$ satisfy $\|f\|_{\infty}\leq 1$. Let $T(x,y) = \inf\{t\geq 0: Y_t(x) = Y_t(y)\} = \inf\{t\geq 0: Y_t(y) - Y_t(x) = 0\}$. By (\ref{eq5.2}) we have
 \beqnn
|Q_tf(x)-Q_tf(y)| \ar\leq\ar \mbf{P}[|f(Y_t(x))-f(Y_t(y))|\mathbf{1}_{\{T(x,y)>t\}}] \cr
 \ar\leq\ar
2\mathbf{P}(T(x,y)\leq t) = 2\mbf{P}[e^{-(y-x)\bar{v}^\xi_{0,t}}].
 \eeqnn
The right-hand side tends to zero as $|x-y|\to 0$. That proves the strong Feller property of $(Q_t)_{t\geq 0}$. \qed

In view of the result of Corollary~\ref{th5.4}, one may naturally expect a characterization of the decay rate of the survival probability $\mbf{P}_x(\tau_0> t)$ as $t\to \infty$. This problem for the CBRE-diffusion was studied by B\"oinghoff and Hutzenthaler (2012). More recently, Bansaye et al.\ (2013) studied the problem for CB-processes with catastrophes, which is actually a CBRE-process with stable branching and finite variation L\'evy environment. Palau and Pardo (2015a) studied the same problem for a CBRE-process with stable branching in a random environment given by a Brownian motion with drift. The decay rate of the survival probability for a CBRE-process with stable branching and a general L\'evy environment was studied in Li and Xu (2016). The strong Feller property of classical CBI-processes was proved in Li and Ma (2015).


\section{CBIRE-processes}

\setcounter{equation}{0}

In this section, we discuss the CBIRE-process defined by (\ref{eq1.3}). Let $h\ge 0$ be a constant and $(1\land u)n(du)$ a finite measure on $(0,\infty)$. Suppose that $(\Omega, \mcr{F},\mcr{F}_t,\mbf{P})$ is a filtered probability space satisfying the usual hypotheses. Let $\{W(t): t\ge 0\}$, $\{L(t): t\ge 0\}$ and $M(ds,dz,du)$ be as before. In addition, let $\{\eta(t): t\ge 0\}$ be an increasing $(\mcr{F}_t)$-L\'evy process with
 \beqlb\label{eq6.1}
\mbf{P} (e^{-\lambda\eta(t)}) = e^{-t\psi(\lambda)}, \qquad \lambda\ge 0,
 \eeqlb
where
 \beqlb\label{eq6.2}
\psi(\lambda) = h\lambda + \int_0^{\infty}(1-e^{-\lambda u})n(du).
 \eeqlb
We assume that all those noises are independent of each other. The construction and basic properties of the CBIRE-process are provided by the following results. We here omit their proofs since they follow by modifications of the arguments in Sections~3 and~4.

\begin{theorem}\label{th6.1} For any positive $\mcr{F}_0$-measurable random variable $Y(0)$, there is a unique positive strong solution $\{Y(t): t\ge 0\}$ to (\ref{eq1.3}).
\end{theorem}

\begin{proposition}\label{th6.2} Let $\{Y(t): t\ge 0\}$ be defined by (\ref{eq1.3}) and $Z(t) = Y(t)\exp\{-\xi(t)\}$ for $t\ge 0$. Then we have
 \beqnn
Z(t) \ar=\ar Y(0) - b\int_0^t e^{-\xi(s)}Y(s) ds + \int_0^t e^{-\xi(s)}\sqrt{2cY(s)} dW(s) \cr
 \ar\ar
+ \int_0^t\int_0^\infty\int_0^{Y(s-)} ze^{-\xi(s-)} \tilde{M}(ds,dz,du) + \int_0^t e^{-\xi(s)} d\eta(s).
 \eeqnn
\end{proposition}

\begin{theorem}\label{th6.3} Let $\mbf{P}^\xi$ be the conditional law given $\{\xi(t): t\ge 0\}$. Then for any $\lambda\geq 0$ and $t\geq r\geq 0$, we have
 \beqnn
\mbf{P}^\xi [e^{-\lambda Y(t)}|\mcr{F}_r]
 =
\exp\bigg\{-Y(r)v^\xi_{r,t}(\lambda) - \int_r^t\psi(v^\xi_{s,t}(\lambda))ds\bigg\}.
 \eeqnn
\end{theorem}

\begin{theorem}\label{th6.4} The Markov process $\{Y(t): t\ge 0\}$ defined by (\ref{eq1.3}) has Feller transition semigroup $(\bar{Q}_t)_{t\ge 0}$ defined by
 \beqlb\label{eq6.3}
\int_{[0,\infty)} e^{-\lambda y} \bar{Q}_t(x,dy)
 =
\mbf{P}\bigg[\exp\bigg\{-xv^\xi_{0,t}(\lambda) - \int_0^t\psi(v^\xi_{s,t}(\lambda))ds\bigg\}\bigg].
 \eeqlb
\end{theorem}

\begin{theorem}\label{th6.5}
Under Grey's condition, the transition semigroup $(\bar{Q}_t)_{t\geq 0}$ defined by (\ref{eq6.3}) has the strong Feller property. \end{theorem}

The transition semigroup $(\bar{Q}_t)_{t\ge 0}$ given by (\ref{eq6.3}) generalizes the one defined by (\ref{eq3.9}). We can give another useful characterization of this semigroup. For this purpose, let us consider an independent copy $\{L'(t): t\ge 0\}$ of the L\'evy process $\{L(t): t\ge 0\}$ with $L'(0)=0$. We may assume the process has the following L\'evy-It\^o decomposition:
 \beqnn
L'(t) = \beta t + \sigma B'(t) + \int_0^t \int_{[-1,1]} (e^z-1)\tilde{N}'(ds,dz) + \int_0^t \int_{[-1,1]^c} (e^z-1) N'(ds,dz),
 \eeqnn
where $\{B'(t): t\ge 0\}$ is a Brownian motion and $N'(ds,dz)$ is a Poisson random measure on $(0,\infty)\times \mbb{R}$ with intensity $ds\nu(dz)$. Let $\{\xi'(t): t\ge 0\}$ be the L\'evy process defined by
 \beqnn
\xi'(t) = at + \sigma B'(t) + \int_0^t \int_{[-1,1]} z\tilde{N}'(ds,dz) + \int_0^t \int_{[-1,1]^c} z N'(ds,dz).
 \eeqnn
Set $L(t) = L(0)-L'(-t-)$ and $\xi(t) = \xi(0)-\xi'(-t-)$ for $t<0$. Then $\{L(t): -\infty< t< \infty\}$ and $\{\xi(t): -\infty< t< \infty\}$ are time homogeneous L\'evy processes. We can extend (\ref{eq2.12}) easily to $r\le t\in \mbb{R}$. In particular, for any $\lambda\ge 0$ there is a unique positive solution $r\mapsto v^\xi_{r,0}(\lambda)$ to
 \beqnn
v^\xi_{r,0}(\lambda) = \lambda-\int_r^0 \phi(v^\xi_{s,0}(\lambda))ds + \int_r^0 v^\xi_{s,0}(\lambda) L(\overleftarrow{ds}), \qquad r\le 0.
 \eeqnn
The result of Theorem~\ref{th2.3} can be extended to $r\le t\in \mbb{R}$. Then $r\mapsto v^\xi_{r,0}(\lambda)$ is also the unique positive solution to
 \beqlb\label{eq6.4}
v^\xi_{r,0}(\lambda) = e^{\xi(0)-\xi(r)}\lambda - \int_r^0 e^{\xi(s)-\xi(r)}\phi(v^\xi_{s,0}(\lambda)) ds, \qquad r\le 0.
 \eeqlb
It follows that $r\mapsto u^\xi_{r,0}(\lambda) := e^{\xi(r)} v^\xi_{r,0}(\lambda)$ is the unique positive solution to
 \beqlb\label{eq6.5}
u^\xi_{r,0}(\lambda) = \lambda - \int_r^0 e^{\xi(s)}\phi(e^{-\xi(s)}u^\xi_{s,t}(\lambda))ds, \qquad r\le 0.
 \eeqlb
By the time homogeneity of the L\'evy process $\{\xi(t): -\infty< t< \infty\}$, we have
 \beqlb\label{eq6.6}
\int_{[0,\infty)} e^{-\lambda y} \bar{Q}_t(x,dy)
 =
\mbf{P}\bigg[\exp\bigg\{-xv^\xi_{-t,0}(\lambda) - \int_{-t}^0\psi(v^\xi_{s,0}(\lambda))ds\bigg\}\bigg].
 \eeqlb
In the subcricital case, a necessary and sufficient condition for the ergodicity of the transition semigroup $(\bar{Q}_t)_{t\geq 0}$ defined by (\ref{eq6.3}) is provided by the following theorem:

\begin{theorem}\label{th6.6} Suppose that $a_1 := \mbf{P}[\xi(1)]< b$. Then there is a probability measure $\mu$ on $[0,\infty)$ so that $\bar{Q}_t(x,\cdot)\to \mu$ weakly as $t\rightarrow \infty$ for every $x\geq 0$ if and only if
 \beqnn
\int_1^{\infty} \log(u)n(du)< \infty.
 \eeqnn
Under the above condition, we have
 \beqlb\label{eq6.7}
\int_{[0,\infty)} e^{-\lambda y} \mu(dy)
 =
\mbf{P}\bigg[\exp\bigg\{-\int_{-\infty}^0\psi(v^\xi_{s,0}(\lambda))ds\bigg\}\bigg].
 \eeqlb
\end{theorem}

\proof Under the assumption, we may adjust the parameters in (\ref{eq1.3}) so that $b> 0> a_1$. Then $u^\xi_{-t,0}(\lambda)\le e^{-bt}\lambda$ by Gronwall's inequality, and hence $u^\xi_{-t,0}(\lambda)\to 0$ as $t\to \infty$. In view of (\ref{eq6.3}), by applying Theorem~1.20 in Li (2011) and dominated convergence we conclude that $\bar{Q}_t(x,\cdot)$ converges to a probability measure $\mu$ as $t\rightarrow \infty$ for every $x\geq 0$ if and only if a.s.
 \beqnn
f(\xi,\lambda) := \int_{-\infty}^0 \psi(v^\xi_{s,0}(\lambda)) ds
 =
\int_{-\infty}^0\psi(e^{-\xi(s)}u^\xi_{s,0}(\lambda))ds< \infty, \qquad \lambda\ge 0.
 \eeqnn
Clearly, $\mu$ is given by (\ref{eq6.7}) if the above condition is satisfied. For any $z>0$, define $\tau_\lambda(z) = \sup\{r<0: u^\xi_{r,0}(\lambda)\leq z\}$. By (\ref{eq6.5}) and a change of the variable, we have
 \beqnn
f(\xi,\lambda) = \int_{-\infty}^0 \frac{\psi(e^{-\xi(s)}u^\xi_{s,0}(\lambda))} {\phi(e^{-\xi(s)}u^\xi_{s,0}(\lambda))} e^{-\xi(s)}du^\xi_{s,0}(\lambda)
 =
\int_0^\lambda \frac{\psi(e^{-\xi(\tau_\lambda(z))}z)} {\phi(e^{-\xi(\tau_\lambda(z))}z)} e^{-\xi(\tau_\lambda(z))}dz.
 \eeqnn
Since $\xi(t)\rightarrow \infty$ as $t\to -\infty$, we have a.s. $M := \sup_{s\le 0} e^{-\xi(s)}< \infty$. It is simple to see that $\phi(z) = bz + o(z)$ as $z\rightarrow 0$. Then $f(\xi,\lambda)< \infty$ if and only if
 \beqlb\label{eq6.8}
\lefteqn{\int_0^\lambda dz\int_0^{\infty}\frac{1-\exp\{-e^{-\xi(\tau_\lambda(z))}zu\}}{z} n(du)}\qquad\ar\ar\cr
 \ar=\ar
\int_0^{\infty} n(du)\int_0^{\lambda u} \frac{1-\exp\{-e^{-\xi(\tau_\lambda(y/u))}y\}}{y}dy< \infty.
 \eeqlb
For all $u>0$ we have $1-\exp\{-e^{-\xi(\tau_\lambda(y/u))}y\}\le My$. It follows that
 \beqnn
\int_{(0,1]} n(du)\int_0^{\lambda u} \frac{1-\exp\{-e^{-\xi(\tau_\lambda(y/u))}y\}}{y}dy
 \le
M\lambda \int_{(0,1]} u n(du)< \infty.
 \eeqnn
For $u>1$ we have $1 - \exp\{-e^{-\xi(\tau_\lambda(y/u))}y\}\to 1$ as $y\to \infty$. Then (\ref{eq6.8}) holds if and only if
 \beqnn
\int_{(1,\infty)} n(du)\int_0^{\lambda u} \frac{1}{y}dy
 =
\int_{(1,\infty)} \log (\lambda u)n(du)< \infty.
 \eeqnn
That implies the desired result. \qed

\begin{corollary}\label{th6.7} Suppose that $a_1 := \mbf{P}[\xi(1)]< b$. Let $(\tilde{Q}_t)_{t\ge0}$ be the transition semigroup defined by (\ref{eq3.9}). Then there is a probability measure $\mu$ on $[0,\infty)$ so that $\tilde{Q}_t(x,\cdot)\to \mu$ weakly as $t\rightarrow \infty$ for every $x\geq 0$ if and only if
 \beqnn
\int_1^{\infty} u\log(u)m(du)< \infty.
 \eeqnn
Under the above condition, we have
 \beqnn
\int_{[0,\infty)} e^{-\lambda y} \mu(dy)
 =
\mbf{P}\bigg[\exp\bigg\{-\int_{-\infty}^0\psi(v^\xi_{s,0}(\lambda))ds\bigg\}\bigg].
 \eeqnn
\end{corollary}

\proof One can see that (\ref{eq3.9}) is the special form of (\ref{eq6.3}) with $\psi = \phi_0'$. Then we get the results by Theorem~\ref{th6.6}. \qed

\end{document}